\documentclass[a4paper]{article}

\usepackage{fullpage}
\usepackage{amsmath,amssymb,amsfonts}
\usepackage{nicefrac}
\usepackage{xcolor}
\definecolor{cb2blue}{RGB}{55,126,184}
\definecolor{cb2green}{RGB}{77,175,74}
\definecolor{cb2red}{RGB}{228,26,28}
\usepackage{subcaption}
\usepackage{graphicx}

\usepackage[ruled,vlined,linesnumbered,algosection,algo2e]{algorithm2e} % Before cleveref
\SetKw{KwAnd}{and}

\usepackage{hyperref}
\hypersetup{colorlinks,
	bookmarks,
	citecolor=cb2green,
	linkcolor=cb2blue,
	urlcolor=cb2red}

\newcommand{\emailLink}[1]{\href{mailto:#1}{#1}}
\newcommand{\orcidLink}[1]{\href{https://orcid.org/#1}{#1}}

\DeclareMathOperator*{\argmin}{\arg\min}
\newcommand{\N}{\mathbb{N}}
\newcommand{\R}{\mathbb{R}}
\newcommand{\coloneqq}{:=}
\newcommand{\identity}{\mathbb{I}}
\DeclareMathOperator{\diag}{diag}
\DeclareMathOperator*{\minimize}{minimize}
\DeclareMathOperator{\stt}{subject~to}
\DeclareMathOperator{\rank}{rank}

\newcommand{\LL}{\mathcal{L}}%Lagrangian
\newcommand{\pdLL}{\mathcal{M}}% primal-dual Lagrangian
\newcommand{\baseqpsolver}{\texttt{QProxIP}}
\newcommand{\rlqpsolver}{\texttt{RL-QProxIP}}

\newcommand\blfootnote[1]{%
	\begingroup
	\renewcommand\thefootnote{}\footnote{#1}%
	\addtocounter{footnote}{-1}%
	\endgroup
}

\newcommand{\TheAuthorJB}{Jeremy Bertoncini}
\newcommand{\TheEmailJB}{jeremy.bertoncini@unibw.de}
\newcommand{\TheOrcidJB}{0009-0003-6093-7532}

\newcommand{\TheAuthorADM}{Alberto De~Marchi}
\newcommand{\TheEmailADM}{alberto.demarchi@unibw.de}
\newcommand{\TheOrcidADM}{0000-0002-3545-6898}

\newcommand{\TheAuthorSG}{Simon Gottschalk}
\newcommand{\TheEmailSG}{simon.gottschalk@unibw.de}
\newcommand{\TheOrcidSG}{0000-0003-4305-5290}

\newcommand{\TheAuthorMG}{Matthias Gerdts}
\newcommand{\TheEmailMG}{matthias.gerdts@unibw.de}
\newcommand{\TheOrcidMG}{0000-0001-8674-5764}

\newcommand{\TheFunding}{This research is funded by dtec.bw -- Digitalization and Technology Research Center of the Bundeswehr [MissionLab, SeRANIS]. dtec.bw is funded by the European Union -- NextGenerationEU}

\newcommand{\TheTitle}{Reinforcement Learning for Adaptive Interior Point Methods in Convex Quadratic Programming}

\newcommand{\TheKeywords}{%
Quadratic programming,
Reinforcement learning,
Amortized optimization,
Solver tuning.%
}

\begin{document}

\title{\bfseries \TheTitle}

\author{%
\TheAuthorJB\thanks{\textsc{email} \emailLink{\TheEmailJB}, \textsc{orcid} \orcidLink{\TheOrcidJB}.}%
\and\TheAuthorADM\thanks{\textsc{email} \emailLink{\TheEmailADM}, \textsc{orcid}  \orcidLink{\TheOrcidADM}.}%
\and\TheAuthorMG\thanks{\textsc{email} \emailLink{\TheEmailMG}, \textsc{orcid}  \orcidLink{\TheOrcidMG}.}%
\and\TheAuthorSG\thanks{\textsc{email} \emailLink{\TheEmailSG}, \textsc{orcid}  \orcidLink{\TheOrcidSG}.}%
}%

\date{}

\maketitle

\blfootnote{%
The authors are with the Department of Aerospace Engineering, Institute of Applied Mathematics and Scientific Computing, University of the Bundeswehr Munich, Werner-Heisenberg-Weg 39, 85577 Neubiberg, Germany.
}%

\blfootnote{\TheFunding.}

\begin{abstract}
    Quadratic programming is a workhorse of modern nonlinear optimization, control, and data science.
    Although regularized methods offer convergence guarantees under minimal assumptions on the problem data, they can exhibit the slow tail-convergence typical of first-order schemes, thus requiring many iterations to achieve high-accuracy solutions.
    Moreover, hyperparameter tuning significantly impacts the solver performance but how to find an appropriate parameter configuration remains an elusive research question.
    To address these issues, we explore how data-driven approaches can accelerate the solution process.
    Aiming at high-accuracy solutions, we focus on a regularized interior-point solver and carefully handle its two-loop flow and control parameters.
    We will show that reinforcement learning can make a significant contribution to facilitating the solver tuning and to speeding up the optimization process.
    Numerical experiments demonstrate that, after a lightweight training, the learned policy generalizes well to different problem classes with varying dimensions.

    \medskip

    \noindent
    \textbf{Keywords}
    \TheKeywords
\end{abstract}

\section{Introduction}

Quadratic programming provides a fundamental optimization model with applications in finance, data analysis, robotics, process control, and operations research.
A quadratic program (QP) with $n$ variables, $m$ equality constraints, and $p$ inequality constraints can be written in the form
\begin{equation}
    \label{eq:QP}
    \minimize_{x \in \R^n}{}\quad \frac{1}{2} x^\top Q x + q^\top x
    \qquad
    \stt{}\quad Ax=b ,\quad Gx \leq d
\end{equation}
where $x\in\R^n$ is the optimization variable, $Q\in\R^{n\times n}$ is a symmetric positive semi-definite matrix that defines the quadratic cost, $q\in\R^n$ is the linear cost vector,
matrix $A\in\R^{m\times n}$ with vector $b\in\R^m$ define the linear equality constraints,
and matrix $G\in\R^{p\times n}$ with vector $d\in\R^p$ represent the linear inequality constraints.

Methods for the numerical solution of QPs have been extensively studied.
These differ in how they balance the number of iterations and the computational cost for each iteration.
First-order methods, such as OSQP \cite{stellato2020osqp} based on ADMM \cite{boyd2011distributed}, build upon simple steps but typically require many iterations, so several acceleration schemes have been proposed \cite{giselsson2017linear,ichnowski2021rlqp,daspremont2021acceleration}.
In contrast, Newton-type methods such as interior point (IP) \cite{gondzio2012interior},
active set \cite{ferreau2014qpoases} and semismooth Newton \cite{liaomcpherson2020fbstab} usually need fewer, yet more demanding, iterations to reach accurate solutions.
Recently, proximal techniques \cite{parikh2014proximal} and the augmented Lagrangian (AL) framework \cite{rockafellar1976augmented} have seen a resurgent interest due to their inherent regularization properties, resulting in the development of solvers such as QPDO \cite{demarchi2022qpdo}, QPALM \cite{hermans2022qpalm} and ProxQP \cite{bambade2025proxqp}.
On this vein, the integration of proximal regularization with the IP strategy led to IP-PMM \cite{pougkakiotis2021interior}, PS-IPM \cite{cipolla2023proximal} and PIQP \cite{schwan2023piqp}.
Thanks to proximal regularization, these solvers can cope with matrix $Q$ merely positive semi-definite and with matrices $A$, $G$ having linearly dependent rows.

Although convergence guarantees can be established under minimal requirements on the problem data and algorithmic parameters,
there are often numerous hyperparameters that must be tuned to speed up convergence, and few clues on how to adjust them.
For example, the barrier parameter in interior-point methods must vanish; since there is no clear indication about its optimal rate of decay, heuristics are typically adopted \cite{schwan2023piqp}.
In proximally regularized methods, the stepsize parameters not only affect the convergence rate but they also impact on the stability of linear algebra operations; cf. \cite[\S 4]{cipolla2023proximal}.
Some theoretical works seek optimal values for the stepsize parameters, but they rely on solving auxiliary problems that are harder than the original QP itself \cite{giselsson2017linear,goujaud2024pepit}.
Another offline tuning approach is that of OPAL \cite{audet2014optimization}, which takes an algorithm as input and returns a recommendation on parameter values that maximize some user-defined performance metric.
Alternatively, online heuristics introduce ``feedback'' by monitoring progress and adapting the parameters along the optimization process;
see \cite[\S 2]{birgin2012augmented}, \cite[\S 5.2]{stellato2020osqp}, \cite[\S 5.3]{hermans2022qpalm} for some common instances.

After the seminal work \cite{li2017learning} on ``learning to optimize'' for automated algorithm design, 
several researchers approached acceleration, warm-starting and hyperparameter tuning as machine learning tasks.
A modern view on this computational approach is presented in the tutorial \cite{amos2023tutorial} on ``amortized optimization'' methods, which use learning to predict the solutions to problems by exploiting the shared structure between similar problem instances.
In recent years, machine learning techniques have been increasingly adopted to improve optimization strategies and to identify new ones \cite{castera2025learning}.
For example, neural networks are used to warm-start QP solvers \cite{sambharya2023end} and fixed-point algorithms \cite{sambharya2024learning}.

Despite the flexibility of learning-based methodologies, considerable attention was devoted to accelerating \emph{first-order} optimization methods, driven by their scalability and applications in machine learning.
Exemplarily, the architecture proposed in \cite{sambharya2023end} consists of a neural network mapping from the QP problem data to warm-starts, learned across families of QPs.
Dedicated to proximal methods for image processing, the tuning-free approach of \cite{wei2022tfpnp} formulates the online hyperparameter selection as a sequential decision-making problem and implements a policy network for automated parameter search.
Closer to our focus here is RLQP \cite{ichnowski2021rlqp}, which couples the OSQP solver with a reinforcement learning (RL) agent:
by learning a policy to adapt the internal parameters of the ADMM algorithm, RLQP can outperform OSQP (with its fine-tuned heuristic update rules) on a variety of problems.

\begin{figure}[tbh]
    \centering%
    \includegraphics{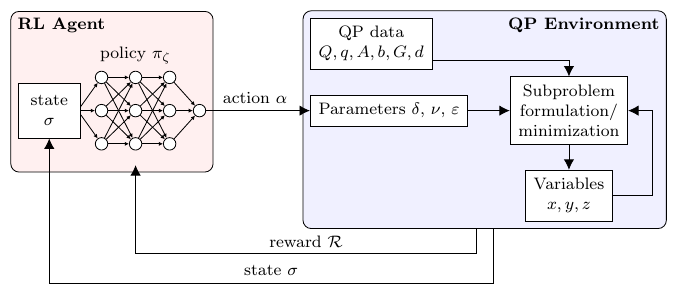}
    \caption{A QP solver is combined with a reinforcement learning (RL) agent that evaluates a policy to adapt the internal parameters of the QP solver.
    Based on its observation of the state of the environment, the RL action aims at improving the solver's performance.
    The RL policy is parameterized by a small artificial neural network, the state $\sigma$ summarizes the internal state of the QP solver (at the beginning of each inner loop), the action $\alpha$ affects the parameters $\delta$, $\nu$, $\varepsilon$ of the solver, and the reward $\mathcal{R}$ promotes faster convergence to a solution of the QP.
    Illustration inspired by \cite[Fig. 1]{ichnowski2021rlqp}.}%
    \label{fig:framework}%
\end{figure}

In contrast, the combination of learning techniques and \emph{second-order}, or Newton-type, methods remains relatively unexplored.
We investigate the coupling of a regularized IP solver for QPs with a policy for online hyperparameter tuning, parametrized with a simple neural network and learned by RL.
Figure~\ref{fig:framework} illustrates the algorithmic framework and the feedback interconnection between QP solver and RL agent.
Following the approach of RLQP \cite{ichnowski2021rlqp}, we learn a policy $\pi_\zeta$ (e.g., a neural network parameterized by $\zeta$) that maps states $\sigma$ to actions $\alpha$ with the goal of maximizing an accumulated reward $\mathcal{R}$.
Compared to the single-loop structure of OSQP, the architecture arising from an IP scheme (the \emph{QP Environment}) can be represented with two nested loops, where the inner loop corresponds to the approximate solution of a subproblem (by executing some iterations of Newton method).
While ubiquitous in second-order methods and with well established convergence theory, the two-loop structure does not fit the RL framework developed for RLQP, because of endogenous parameter sequences controlling the solution process (and affecting the MDP formalization).
We aim to fill this gap by exploring a machine learning framework suitable for accelerating the tuning and solution process of QP solvers beyond first-order methods.

\paragraph{Contribution}

We formulate a RL framework to learn a policy for effective and reliable hyperparameter tuning, with the base QP solver retaining provable convergence guarantees.
The definitions of state, action and reward are intended (i) to yield a policy that is invariant to problem size, scaling and permutations,
and (ii) to capture the two-loop structure typical of second-order methods.

Experimental results on random and benchmark QPs demonstrate that the learned policy outperforms the base solver and is competitive with state-of-the-art solvers.
Moreover, it generalizes well to previously unseen problem sizes, scales, and classes, where tedious hand-tuning is typically required.
Especially on ill-conditioned problems, we observe that the RL-enhanced solver automatically adapts its behavior, outperforming other solvers in terms of robustness.
In addition, these capabilities are achieved after lightweight training of a small neural network on a standard laptop, in contrast to the substantially higher computational cost associated with training RLQP \cite{ichnowski2021rlqp}.

\paragraph{Outline}
Convex quadratic programming is introduced in Section~\ref{sec:solving_qp} with a stabilized IP numerical scheme and a discussion on the challenges and opportunities of hyperparameter tuning.
An overview of reinforcement learning is given in Section~\ref{sec:rl}, with the tools needed to properly interact with the QP solver.
Training procedures and numerical results for three evaluation settings are then presented in Section~\ref{sec:num_results}.
We conclude with some final remarks in Section~\ref{sec:conclusions}.

\medskip

Throughout this article we use various common acronyms, introduced in the table below.
Our numerical experiments compare a QP solver named \baseqpsolver{}, developed by the authors based on the methodology described in Section~\ref{sec:solving_qp}, and its RL–enhanced extension \rlqpsolver{}. 

\begin{table}[h]
    \centering%
    \begin{tabular}{ll|ll}
        \hline
        LP/QP & Linear/Quadratic Problem & RL & Reinforcement Learning \\
        IP & Interior Point & MDP & Markov Decision Process \\
        AL & Augmented Lagrangian & PPO & Proximal Policy Optimization\\
        PPA & Proximal Point Algorithm & \baseqpsolver{} & base QP solver \\
        ADMM & Alternating Direction Method of Multipliers & \rlqpsolver{}   & RL–enhanced QP solver\\
        \hline
    \end{tabular}
\end{table}

\section{Solving convex QPs}\label{sec:solving_qp}

In view of IP methods, the QP problem \eqref{eq:QP} can be equivalently formulated using an auxiliary variable $s\in\R^p$ as 
\begin{equation}\label{eq:QP_with_s}
    \minimize_{x \in \R^n, s \in \R^p}\quad
    \frac{1}{2} x^\top Q x + q^\top x
    \qquad
    \stt\quad
    Ax=b ,\quad
    Gx - d +s =0 ,\quad
    s \geq 0 .
\end{equation}
It is a standard result that, for convex QPs, the following necessary conditions are also sufficient for optimality:
there exists Lagrange multipliers $y\in\R^m$ and $z\in\R^p$ such that
\begin{align}\label{eq:optimality_conditions}
    0 ={}& Qx+q + A^\top y + G^\top z ,& 
    0 ={}& \min\{ s, z \} ,&
    0 ={}& Ax-b ,&
    \text{and}\quad
    0 ={}& Gx-d+s ,
\end{align}
where the second equality coincides with the complementarity condition $0 \leq s \perp z \geq 0$.
Therefore, a valid criterion for returning a triplet $(x,y,z)$ as approximately optimal, with tolerance $\epsilon \geq 0$,
is to check that the conditions
\begin{align}\label{eq:termination_conditions}
    r_{\rm prim}(x,y,z) \coloneqq{}& \left\| \begin{array}{c}
        Ax-b \\
        \min\{ d-Gx, z \}
    \end{array} \right\|_\infty \leq \epsilon &
    \text{and}\quad
    r_{\rm dual}(x,y,z) \coloneqq{}& \left\| Qx + q + A^\top y + G^\top z \right\|_\infty \leq \epsilon
\end{align}
are satisfied.
The primal $r_{\rm prim}$ and dual $r_{\rm dual}$ residual functions indicate how far a triplet $(x,y,z)$ is from being a primal-dual solution to \eqref{eq:QP}, in view of \eqref{eq:optimality_conditions}.
As such, these metrics will play a crucial role in monitoring and quantifying the performance of a solution process.

\subsection{Regularized interior point schemes}\label{sec:ip}

Moving towards numerical schemes, we notice that, from an operator perspective, \eqref{eq:optimality_conditions} is equivalent to a monotone inclusion.
Adopting the (inexact) proximal point algorithm (PPA) \cite{rockafellar1976monotone,luque1984asymptotic}, different reformulations of \eqref{eq:optimality_conditions} lead to different regularized schemes,
whose subproblems correspond to the (approximate) evaluation of the proximal operator.
This construction naturally gives rise to two-loop algorithms: the outer for PPA and the inner for the subproblem solution.
Among others, QPDO \cite{demarchi2022qpdo} and PS-IPM \cite{cipolla2023proximal} have been developed in this spirit, so that proximally regularized subproblems can be (efficiently and reliably) solved with Newton-type methods.

Here we outline a regularized IP scheme that, reminiscent of QPDO and PS-IPM, treats equality and inequality constraints respectively with primal-dual augmented Lagrangian and barrier strategies.
For further regularization, a Tikhonov-like proximal term is included to handle also non-strictly convex QPs, as in \cite{bambade2025proxqp,bemporad2018numerically,stellato2020osqp}.

\paragraph{Primal IP}
Beginning from the classical augmented Lagrangian function for \eqref{eq:QP_with_s}, for any given primal-dual estimate $\theta = (\theta_x,\theta_y,\theta_z)\in\R^{n+m+p}$, regularization parameter $\delta=(\delta_x,\delta_y,\delta_z)\in\R_{++}^3$ and barrier parameter $\nu>0$,
we define the merit function $\LL_\delta(\cdot,\cdot;\theta,\nu)\colon \R^n\times\R_{++}^p \to \R$ as
\begin{multline}\label{eq:primal_aug_lagrangian}
    \LL_\delta(x,s;\theta,\nu)
    \coloneqq
    \frac{1}{2} x^\top Q x + q^\top x
    + \theta_y^\top (Ax-b)
    + \theta_z^\top (Gx-d+s) \\
    + \frac{\delta_x}{2} \|x-\theta_x\|^2
    + \frac{1}{2\delta_y} \| Ax-b \|^2
    + \frac{1}{2\delta_z} \| Gx-d+s \|^2
    - \nu \sum_{i=1}^p \ln s_i
    .
\end{multline}
As usual in barrier methods, the positivity constraint $s \geq 0$ in \eqref{eq:QP_with_s} is incorporated in a smooth way by replacing its indicator with a (logarithmic) barrier function: the smaller the barrier parameter $\nu$, the better the approximation.
Note that, for any fixed $\theta$, $\delta$ and $\nu$, $\LL_\delta(\cdot,\cdot;\theta,\nu)$ is strictly convex under mere convexity of \eqref{eq:QP}.
Hence, it always admits a unique minimizer $(x,s)$ that necessarily satisfies $s>0$.
Effectively, the minimization of $\LL_\delta(\cdot,\cdot;\theta,\nu)$ constitutes a suitable subproblem for tackling \eqref{eq:QP}, since its minimizer is closely related to the proximal point of an operator associated to \eqref{eq:optimality_conditions}; see \cite{rockafellar1976augmented,cipolla2023proximal}.

In practice, by minimizing $\LL_{\delta^k}(\cdot,\cdot;\theta^k,\nu_k)$ one can obtain an improved primal estimate $(x^k,s^k)$ and classical multiplier updates from AL methods yield a dual estimate $(y^k,z^k)$ as
\begin{equation}\label{eq:dual_estimate}
    y^k \gets \theta_y^k + \frac{Ax^k-b}{\delta_{y,k}}
    \quad\text{and}\quad
    z^k \gets \theta_z^k + \frac{Gx^k-d+s^k}{\delta_{z,k}}
    .
\end{equation}
After updating the incumbent estimate $\theta^{k+1}$ and the hyperparameters $\delta^{k+1}$ and $\nu_{k+1}$, typically based on $(x^k,s^k,y^k,z^k)$, one moves on to the next subproblem.

\paragraph{Primal-dual IP}
A widespread technique is that of \emph{primal-dual} interior point methods, where the subproblem is modified to involve both primal and dual variables,
so that improved values for $(y,z)$ are directly obtained from the inner minimization procedure.
One way to generate a suitable subproblem for this task is to adopt the primal-dual augmented Lagrangian function \cite{gill2012primal}.
We define $\pdLL_\delta(\cdot,\cdot,\cdot,\cdot;\theta,\nu)\colon\R^n\times\R_{++}^p\times\R^m\times\R^p\to\R$ as
\begin{equation}\label{eq:primal_dual_aug_lagrangian}
    \pdLL_\delta(x,s,y,z;\theta,\nu)
    \coloneqq
    \LL_{\delta}(x,s;\theta,\nu )
    +
    \frac{1}{2 \delta_y} \| Ax-b + \delta_y (\theta_y - y) \|^2 + \frac{1}{2 \delta_z} \| Gx-d+s + \delta_z (\theta_z - z) \|^2
\end{equation}
where the additional terms
guarantee that, when minimizing $\pdLL_{\delta^k}(\cdot,\cdot,\cdot,\cdot;\theta^k,\nu_k)$ exactly, the optimal values for $(y,z)$ coincide with the classical update \eqref{eq:dual_estimate} that would be performed within a purely primal scheme.
Similarly to its primal counterpart, the primal-dual augmented Lagrangian function $\mathcal{M}$ depends on the regularization parameters $\delta$ and barrier parameter $\nu$.

\begin{algorithm2e}[tbh]
    \DontPrintSemicolon
    \KwData{primal-dual guess $(x,y,z)$, tolerance $\epsilon$}
    set $\theta^0 \gets (x,y,z)$\;
    \For{$k = 0,1,2\ldots$}{
        choose $\delta^k, \nu_k, \varepsilon_k>0$\tcp*{hyperparameters update}\label{step:rpdip:hyperparameters}
        define $\pdLL_k\coloneqq (x,s,y,z) \mapsto \pdLL_{\delta^k}(x,s,y,z;\theta^k,\nu_k)$\;
        find $(x^k,s^k,y^k,z^k)$ such that $\left\Vert \nabla \pdLL_k(x^k,s^k,y^k,z^k) \right\Vert \leq \varepsilon_k$ and $s^k>0$\tcp*{subproblem solve}\label{step:rpdip:subproblem}
        \If(\tcp*[f]{termination check}){$r_{\rm prim}(x^k,y^k,z^k) \leq \epsilon$ \KwAnd $r_{\rm dual}(x^k,y^k,z^k) \leq \epsilon$\label{step:rpdip:termination}}{\KwRet $(x^k,y^k,z^k)$}
        set $\theta^{k+1} \gets (x^k,y^k,z^k)$\tcp*{estimates update}\label{step:rpdip:theta}
    }
    \caption{Regularized interior point method for convex QPs.}
    \label{alg:rpdip}
\end{algorithm2e}

The resulting primal-dual scheme is outlined in Algorithm~\ref{alg:rpdip}, where Step~\ref{step:rpdip:subproblem} poses the most computational requirements and will be discussed in the following Section~\ref{sec:inner_newton}.
The termination conditions \eqref{eq:termination_conditions} are checked at Step~\ref{step:rpdip:termination}, before updating the incumbent primal-dual estimate $\theta^{k+1}$ at Step~\ref{step:rpdip:theta}.
Of all steps, the choice of hyperparameters $\delta^k$, $\nu_k$ and $\varepsilon_k$ at Step~\ref{step:rpdip:hyperparameters} has the least guidelines or specifications dictated by the theory.
Their tuning and influence on the solution process are discussed in Section~\ref{sec:hyperparameters} below.

Convergence of Algorithm~\ref{alg:rpdip} is guaranteed insofar as the regularization parameters $\delta^k$ remain positive and bounded,
the barrier parameter $\nu_k$ vanishes from above,
and the inner tolerance $\varepsilon_k$ is summable.
The technical requirements and mild conditions for convergence are detailed, e.g., in \cite{luque1984asymptotic}.
The convergence properties of Algorithm~\ref{alg:rpdip} can be analyzed with the lens of inexact proximal point algorithms;
we refer the interested reader to \cite{cipolla2023proximal,pougkakiotis2021interior,demarchi2023regularized,demarchi2022qpdo,liaomcpherson2020fbstab,hermans2022qpalm,bemporad2018numerically,bambade2025proxqp} for recent results established by exploiting this connection for regularized QP solvers.

\subsection{Inner minimization with Newton}\label{sec:inner_newton}

The minimization procedure adopted to solve the subproblem at Step~\ref{step:rpdip:subproblem} of Algorithm~\ref{alg:rpdip} is typically based on Newton's method with a backtracking line search procedure for globalization.
Methods of this kind for convex unconstrained minimization have been extensively studied and provide strong convergence guarantees.
However, a direct application of Newton method to minimize $\pdLL_k$ is usually deemed impractical and avoided, since the Hessian matrix $\nabla^2 \pdLL_k$ tends to be dense, but, owing to convexity of $\pdLL_k$, it is sufficient to find a zero of $\nabla \pdLL_k$.
Using the linear transformation $T_{\delta^k}$ as in \cite[\S 3.1]{demarchi2022qpdo}, we can obtain an equivalent root-finding subproblem that retains structure and sparsity of the original problem \eqref{eq:QP}:
\begin{equation}
    0
    =
    T_{\delta^k} \nabla \pdLL_k(x,s,y,z)
    =
    \begin{pmatrix}
         Q x + q + A^\top y + G^\top z + \delta_{x,k} (x-\theta_x^k) \\
         - \nu_k s^{-1} + z \\
         Ax - b + \delta_{y,k} (\theta_y^k - y) \\
         Gx-d+s + \delta_{z,k} (\theta_z^k - z)
    \end{pmatrix}
    ,
\end{equation}
where the vector-inverse in the second block-row is intended componentwise.
Scaling the second block by $s > 0$, the classical primal-dual IP subproblem is recovered with additional proximal regularization terms, as in \cite[\S 2]{pougkakiotis2021interior}, \cite[\S 3]{cipolla2023proximal}.
Then, at any point $(x,s,y,z)$, a Newton direction $(\Delta x,\Delta s,\Delta y,\Delta z)\in\R^{n+p+m+p}$ can be obtained by solving the following system of linear equations:
\begin{equation}\label{eq:newton_system}
    \begin{bmatrix}
        Q + \delta_{x,k} \identity_n & 0 & A^\top & G^\top \\
        0 & \diag(z) & 0 & \diag(s) \\
        A & 0 & -\delta_{y,k} \identity_m & 0 \\
        G & \identity_p & 0 & -\delta_{z,k} \identity_p
    \end{bmatrix}
    \begin{pmatrix}
        \Delta x \\
        \Delta s \\
        \Delta y \\
        \Delta z
    \end{pmatrix}
    =
    -
    \begin{pmatrix}
         Q x + q + A^\top y + G^\top z + \delta_{x,k} (x-\theta_x^k) \\
         - \nu_k + s \odot z \\
         Ax - b + \delta_{y,k} (\theta_y^k - y) \\
         Gx-d+s + \delta_{z,k} (\theta_z^k - z)
    \end{pmatrix}
    ,
\end{equation}
for which there exist several numerical techniques.
In \eqref{eq:newton_system}, $\identity_n$ denotes the identity matrix of size $n$, $\diag(s)$ builds a square diagonal matrix with diagonal $s$, and $\odot$ indicate the Hadamard product.
As our focus is on the combination of the optimization routine in Algorithm~\ref{alg:rpdip} with a learning-based hyperparameter policy,
the (sparse) linear system \eqref{eq:newton_system} is solved directly through LU decomposition with pivoting, for numerical stability, without any additional structure exploitation.
It is worth mentioning here that the regularization parameters $\delta$ are introduced as penalty weights in the definition \eqref{eq:primal_aug_lagrangian} of the augmented Lagrangian,
with the goal of discouraging constraint violations,
but they also play an important role in mitigating the ill-conditioning of matrices arising from the IP scheme; cf. \cite[\S 4]{cipolla2023proximal}.

\subsection{Hyperparameters}\label{sec:hyperparameters}

Although convergence guarantees are retained under minimal assumptions, the practical behavior and performance of Algorithm~\ref{alg:rpdip} depend on the sequences of hyperparameters $\{\delta^k\}$, $\{\nu_k\}$ and $\{\varepsilon_k\}$.
While parameters $\nu$ and $\varepsilon$ should become sufficiently small (relative to a user-specified tolerance), there is no clear way to choose the weights $\delta$;
in fact, convergence is guaranteed even if they are kept constant.

As depicted in Figure~\ref{fig:framework}, all hyperparameters will be subject to the action of the external RL agent and updated according to its policy.
In particular, at Step~\ref{step:rpdip:hyperparameters} of Algorithm~\ref{alg:rpdip}, the regularization parameters $\delta$, the barrier parameter $\nu$, and the inner tolerance $\varepsilon$ are all updated through the output $\alpha$ of the trained RL policy.
Our approach is to define the update rule as
\begin{equation}\label{eq:delta_update}
\begin{aligned}
    \delta_{x,k+1} ={}& \max\{ \delta_{\min}, \alpha_{x,k} \delta_{x,k} \}, &
    \delta_{y,k+1} ={}& \max\{ \delta_{\min}, \alpha_{y,k} \delta_{y,k} \} , &
    \delta_{z,k+1} ={}& \max\{ \delta_{\min}, \alpha_{z,k} \delta_{z,k} \} ,\\
    \nu_{k+1} ={}& \max\{ \nu_{\min}, \alpha_{\nu,k} \nu_{k} \}, &
    \varepsilon_{k+1} ={}& \max\{ \varepsilon_{\min}, \alpha_{\varepsilon,k} \varepsilon_{k} \}, 
\end{aligned}
\end{equation}
starting from some $\delta^0 = (\delta_{x,0},\delta_{y,0},\delta_{z,0}) \in \R_{++}^3$, $\nu_0>0$ and $\varepsilon_0>0$.
For notational convenience, the single actions $\alpha_{x,k}$, $\alpha_{y,k}$, $\alpha_{z,k}$, $\alpha_{\nu,k}$, $\alpha_{\varepsilon,k}\in(0,1)$ are collected in a vector $\alpha^k$.
The design choice in \eqref{eq:delta_update} of a saturated linear decrease is common; it follows the update rules in QPALM \cite{hermans2022qpalm}, ProxQP \cite{bambade2025proxqp} and PIQP \cite{schwan2023piqp}, to name a few.
The small thresholds $\delta_{\min},\nu_{\min},\varepsilon_{\min}>0$ safeguard the numerical linear algebra when addressing the linear system \eqref{eq:newton_system}.

To fix ideas, we let the RL action take the form
\begin{equation}\label{eq:alpha_update}
    \alpha^k
    =
    \alpha
    \left(
    \bar{r}_{\rm prim}, \bar{r}_{\rm dual}, \nu_k, \varepsilon_k
    \right)
    ,
\end{equation}
where $\bar{r}_{\rm prim}, \bar{r}_{\rm dual}$ denote (normalized) residuals at $(x^k,y^k,z^k)$, and the mapping $\alpha$ is learned by training a neural network $\pi_\zeta$ (that is, optimizing its weight and biases $\zeta$).
The choice of defining this mapping with residuals $\bar{r}_{\rm prim}, \bar{r}_{\rm dual}\in\R_+$ as inputs, as opposed to primal-dual decision variables $(x,y,z)\in\R^{n+m+p}$, is intended to make the mapping $\alpha$ invariant to the problem dimensions $(n,m,p)$; more details will be given in Section~\ref{sec:state_reward}.

We conclude this overview of the QP solver by noting that various features (for numerical stability and practical performance) are implemented in both the base and the RL-enhanced versions of the solver.
These standard techniques include
presolving and Ruiz-based preconditioning \cite{stellato2020osqp},
nonmonotone line search \cite{Grippo186_linesearch},
boundary control \cite{pougkakiotis2021interior},
centrality correction \cite{schwan2023piqp}.

\section{Learning-based parameter adaptation}\label{sec:rl}

In the following, we aim at finding a good strategy to choose the parameters during the optimization.
So far, the choice of these parameters is up to the operators, and thus, the optimization highly depends on their experiences.
Reinforcement learning \cite{sutton1998rl} makes it possible to automate this choice and avoid manual hyperparameter tuning.
Therefore, reinforcement learning approaches benefit from their ability to train a policy by interacting with the environment in a trial-and-error manner.

Compared to previous work \cite{ichnowski2021rlqp} on the RL-acceleration of an ADMM-based QP solver, we focus on the regularized interior point approach discussed earlier. Thus, the goal is to automate the choice of hyperparameters $\delta^k$, $\nu_k$, $\varepsilon_k$ at each iteration $k$, by acting through $\alpha^k$ according to \eqref{eq:delta_update}. Thereby, the choice should be made with respect to scale-invariant properties of the problems.
In our case we choose the primal and dual residuals, the barrier parameter and the inner tolerance as in \eqref{eq:alpha_update}.
This design choice is supported by the numerical experiments presented in Section~\ref{sec:num_results}.
Nevertheless, the residuals naturally give a succinct representation of the optimization status, while the barrier parameter reflects the current relaxation of the problem and the inner tolerance determines the effort required to force convergence. Both directly influence the synchronization between the proximal regularization and the performance of the Newton solves.

Figure~\ref{fig:framework} shows the embedding of the RL policy into the optimizer which is summarized in Algorithm~\ref{alg:rpdip}.
In an inner loop, the augmented Lagrangian function $\pdLL_\delta(\cdot,\cdot,\cdot,\cdot;\theta,\nu)$ is minimized for fixed $\theta=\theta^k$, barrier parameter $\nu=\nu_k$ and weights $\delta=\delta_k$, up to a tolerance $\varepsilon_k$.
Once this has been achieved, the estimates $\theta^k$ are updated for the next outer loop iteration.
In addition, the hyperparameters $\delta^k$, $\nu_k$, and $\varepsilon_k$ are adapted by the factors $\alpha^k$ determined by the RL policy.
Thus, the RL policy acts at Step~\ref{step:rpdip:hyperparameters} in Algorithm~\ref{alg:rpdip}, updating all hyperparameters.
This procedure is repeated until the primal and dual residuals are smaller than a tolerance $\epsilon>0$.
During training, the policy to be learned is improved step by step by interacting with the described framework, which in RL context is also called ``environment'',  while trying to solve some training QPs.

\subsection{State and reward}\label{sec:state_reward}

We continue with the mathematical background, which is needed for RL.
It is based on the Markov decision process (MDP) \cite{feinberg2002handbook}, which is represented by the quadruple $(\mathcal{S},\mathcal{A},\mathcal{P},\mathcal{R})$.
The state space $\mathcal{S}$ comprises all possible configurations of the environment, while the action space $\mathcal{A}$ represents all actions which allow the policy to influence the environment.
The probability transition $\mathcal{P}$ describes the probability of reaching a state under the condition that the previous state and action are given.
Finally, the reward function $\mathcal{R}:\mathcal{S}\times \mathcal{S}\rightarrow \R$ rates the state transition and enables to compare and improve solutions.
While the state and action space as well as the reward function need to be defined by the operator, the probability transition does not need to be given explicitly.
It is compensated by the interaction of the policy with the environment.
In the following, we will define this framework for our present task. 

\paragraph{State}

One reason to use RL in this framework is that it is not clear what is the best choice for the factors $\alpha^k$ and what their choice depends on.
Thus, it is natural to choose a state, which best represents the optimization process.
We are convinced that the primal and dual residuals of the QPs are good representatives of the solver status and decisive factors for choosing $\alpha^k$.
Thus, we make them part of the state definition and interpret the upcoming numerical results as indication supporting our choice.

First, the normalized residuals $\bar{r}_{\rm prim}$ and $\bar{r}_{\rm dual}$ defined by
\begin{subequations}
    \begin{align}
    \bar{r}_{\rm prim}(x,y,z)
    \coloneqq{}&
    \frac{r_{\rm prim}(x,y,z)}{\max\{\|A\|_\infty, \|b\|_\infty, \|G\|_\infty, \|d\|_\infty\}} \\
    \bar{r}_{\rm dual}(x,y,z)
    \coloneqq{}&
    \frac{r_{\rm dual}(x,y,z)}{\max\{\|Q\|_\infty, \|q\|_\infty, \|A\|_\infty, \|G\|_\infty\}}
\end{align}
\end{subequations}
are intended to make the state definition invariant to problem scaling.
Then, since the (normalized) residuals and hyperparameters can span several orders of magnitude during the solution process, causing scaling issues while training the RL policy $\pi_\zeta$, we encode the environment state $\sigma$ with a logarithmic scale:
\begin{equation}
    \sigma(\bar{r}_{\rm prim},\bar{r}_{\rm dual},\nu,\varepsilon)
    \coloneqq
    \begin{pmatrix}
        \sigma_{\rm prim} \\
        \sigma_{\rm dual} \\
        \sigma_\nu\\
        \sigma_\varepsilon
    \end{pmatrix}
    \coloneqq
    \begin{pmatrix}
        - \log \left(\bar{r}_{\rm prim}+10^{-9}\right) \\
        - \log \left(\bar{r}_{\rm dual}+10^{-9}\right) \\
        -\log \left(\nu+10^{-17}\right) \\
        -\log \left(\varepsilon+10^{-9}\right)
    \end{pmatrix} 
     ,
\end{equation}
where the small shifts are included to prevent numerical instabilities during the training process.

\paragraph{Action}
The action of our task is the choice of $\alpha^k$ in each iteration step $k$.
Since $\alpha^k$ contains factors that are used to sequentially decrease the hyperparameters $\delta^k$, $\nu_k$, and $\varepsilon_k$, according to \eqref{eq:delta_update}, it must be positive and smaller than one.

\paragraph{Reward}
The reward $\mathcal{R}$ is constructed as a combination of a continuous term and a single success incentive, namely $\mathcal{R} \coloneqq \mathcal{R}_{\rm base} + \mathcal{R}_{\rm bonus}$ where
\begin{align*}
    \mathcal{R}_{\rm base}
    \coloneqq{}&
    -\frac{1}{2}\left[
    \log(\bar{r}_{\rm prim} + 10^{-9})
    +\log(\bar{r}_{\rm dual} + 10^{-9})
    \right]
    -\frac{N_{\rm in}}{N_{\rm in}^{\max}},
    \\
    \mathcal{R}_{\rm bonus}
    \coloneqq{}&
    \left[
    1 + 4\left(1 - \frac{N_{\rm out}}{N_{\rm out}^{\max}}\right)
    \right]
    \iota(\bar{r}_{\rm prim} \leq \epsilon) \iota(\bar{r}_{\rm dual}\leq \epsilon),
\end{align*}
and the characteristic function $\iota$ is defined as
$\iota(\texttt{cond})=1$ if condition \texttt{cond} is true, $\iota(\texttt{cond})=0$ otherwise.
The base reward $\mathcal{R}{_\text{base}}$ encourages convergence by reduction of residuals while penalizing an excessive number of inner Newton iterations $N_{\text{in}}$.
The bonus term $\mathcal{R}{_\text{bonus}}$ is activated only upon satisfaction of the termination criterion at Step~\ref{step:rpdip:termination} of Algorithm~\ref{alg:rpdip}.
When satisfied, it provides a bounded incentive depending on the number of outer iterations required $N_{\text{out}}$. The resulting formulation remains intentionally simple conforming with the MDP framework.

Note that $N_{\rm in}$ and $N_{\rm out}$ are not explicitly part of the state $\sigma$.
However, under the hypothesis that the state $\sigma$ provides a good representation of the current status of the optimizer, then one can expect similar values of $N_{\rm in}$ and $N_{\rm out}$ for similar states $\sigma$.
This observation leads us to consider $N_{\rm in}$ and $N_{\rm out}$ as (indirectly) dependent on $\sigma$, while relying on the fact that any deviation can be seen as variance of our transition probability.

\paragraph{Neural network}

After depicting the MDP, all that remains is to define the policy in order to have a RL framework.
We parametrize our policy $\pi_{\zeta}$ by means of a feedforward neural network, whose weights and biases form the parameters $\zeta$.
We determined the size of the network by a grid search and decided on a net with two hidden layers with $25$ neurons in each layer.
The size of the input layer (four) and output layer (three) are given by the dimension of the state and action space.
For each layer, the commonly used hyperbolic tangent is chosen as the activation function, because of its zero centering and smoothness.

Based on these definitions, RL aims to maximize the objective with respect to the weights and biases of the neural network:
\begin{equation*}
    \min_{\zeta} \quad \mathbb{E}\left[\sum_{i=0}^{\infty} \gamma^i \mathcal{R}(\sigma_i,\sigma_{i+1})\right]
\end{equation*}
with discount factor $\gamma \in (0,1)$ and $\sigma_i \in \mathcal{S}$ for $i\in \N$.
Note that the expected value, which appears in the objective function, is defined over all trajectories $(\sigma_0,\alpha_0,\sigma_1,\alpha_1,\dots)$ for $\alpha_i \in \mathcal{A}$.

\subsection{Training procedures}
Based on the MDP, which was defined in the previous subsection, we run the RL approach. We chose Proximal Policy Optimization (PPO) \cite{schulman2017proximal}, as it is considered robust and stable and is capable of handling continuous action spaces. It is an actor-critic RL approach, which means that the policy $\pi_{\zeta}$ and the value function $V_{\Phi}$ with parameters $\Phi$ are trained simultaneously. Note that in our case, both the policy and the value function share the majority of the neural network, which was described before. 

For an update step in the training, the policy is used for various training problems within the above introduced framework (see Figure~\ref{fig:framework}).
At the same time, the states and actions are stored.
In order to update the policy, PPO forces that the fraction of the updated policy over the previous policy increases for well-rated and decreases for poorly rated behavior in the stored data.
In this way, good behavior is reinforced. The rating of the behavior is given by the general advantages estimation and the policy update is clipped, which avoids big update steps, for a more robust improvement.
Afterwards, the general advantages estimation is updated with respect to the generated data and its corresponding rewards.

Our technical implementation is based on the RL library \texttt{RLlib} \cite{liang2018rllib}, available at \url{http://rllib.io}.
The choice of hyperparameters for the described PPO approach can be seen in Table \ref{tab:PPO_hyperpara}.
Thereby, the learning rate schedule is tailored to the magnitude of the reward and decreases for the final stages of the training.
The gradient clip and the Kullback-Leibler parameters (KL coefficient, KL target), which together with the clipping parameter $\epsilon_{\text{clip}}$ limit the policy updates, are chosen empirically in order to avoid an instable training.
The choice of the entropy schedule leads together with a high number of samples used for an update step (batch size) to an exploratory, but stable training.
Typical default values have been selected for other PPO parameters, which we have confirmed by the robust training behavior.
These include, for instance, the constant $\lambda$ depicting how strongly future rewards influence the update of the general advantages estimation.

\begin{table}[tbh]
    \centering%
    \begin{tabular}{llll}
        \hline
        Learning rate schedule & 
        \begin{tabular}[c]{@{}l@{}}
            (0, $5\cdot 10^{-4}$) \\
            (500k, $3\cdot 10^{-4}$) \\
            (1M, $10^{-4}$) 
        \end{tabular} 
        & Entropy coeff. schedule &
        \begin{tabular}[c]{@{}l@{}}
            (0, $5\cdot 10^{-3}$) \\
            (1M, $2\cdot 10^{-3}$) 
        \end{tabular} \\
        \hline
        Batch size & $8192$ & $\gamma$ & $0.99$ (default) \\
        $\lambda$ & $1.0$ (default)& $\epsilon_{\text{clip}}$ & $0.2$ (default) \\
        Gradient clip & $25$ & KL coeff. & $0.05$ \\
        KL target & $0.01$ & Network layers & [25, 25] (tanh) \\
        \hline
        System configuration & \multicolumn{3}{l}{Kubuntu 22.04.5 LTS, Intel Core i9-12950HX} \\
        Num. CPUs (local worker) & $20$  & Num. CPUs per worker & $1$ \\
        \hline
    \end{tabular}%
    \caption{PPO hyperparameters used in training.}%
    \label{tab:PPO_hyperpara}%
\end{table}

For the training and numerical experiments reported in this work, the QP solver is configured as summarized in Table~\ref{tab:training_solversetup}.
The tolerance $\epsilon$ in the termination conditions~\eqref{eq:optimality_conditions} is set to $10^{-6}$. The minimum values of $\delta$, $\nu$, $\varepsilon$, and the step-size parameter $\alpha$ define the operational bounds of the algorithm, i.e., its standard numerical regime.

The inner-loop tolerance $\varepsilon_k$ is bounded below by $\epsilon/10$, since the inner subproblems do not need to converge beyond this level to guarantee convergence of the outer iterations.
The parameter $\nu_{\min}$ is chosen close to machine precision to handle near-degenerate cases with almost linear dependence in the constraints, where a very small barrier and a fine relaxation of feasibility are required to obtain accurate solutions.
The regularization lower bound $\delta_{\min} = 10^{-12}$ is applied uniformly to AL methods and proximal terms.
These weights push Newton-type methods toward their numerical limits, while improving performance by strongly penalizing the merit function terms~\eqref{eq:primal_aug_lagrangian}, yielding KKT systems \eqref{eq:newton_system} that remain very close to those of the original problem.
The initialization scales $\delta_{x,0}$, $\delta_{y,0}$, and $\delta_{z,0}$ are tuned manually based on the baseline \baseqpsolver{}.
The initial inner tolerance $\varepsilon_0$ is selected such that it is smaller than the initial residual norm, preventing premature satisfaction of the first outer iteration.
The bounds $\alpha_{\min}$ and $\alpha_{\max}$ are set strictly within $(0,1]$ to ensure monotone decrease.
Finally, the iteration limits $N_{\rm in}^{\max}$ and $N_{\rm out}^{\max}$ constrain the algorithmic budget for each problem instance.

{\begin{table}[tbh]
	\centering%
	\begin{tabular}{ccccccc}
		\hline
		$\epsilon$ & $\delta_{\min}$ & $\nu_{\min}$ & $\varepsilon_{\min}$ & $\delta_{x,0}$ & $\delta_{y,0}$ & $\delta_{z,0}$ \\
		$10^{-6}$ & $10^{-12}$ & $10^{-16}$ & $ 0.1\cdot \epsilon$ & $1$ & $10$ & $10$  \\
		\hline 
	    $\nu_0$ & $\varepsilon_0$ & $\alpha_{\min}$ & $\alpha_{\max}$ & $N_{\rm in}^{\max}$ & $N_{\rm out}^{\max}$ &  \\
		$1$ & $ 0.1\cdot \min\{r_{\mathrm{prim}}(\theta^0), r_{\mathrm{dual}}(\theta^0)\}$ & $0.05$ & $0.95$ & $25$ & $25$ &  \\
		\hline
	\end{tabular}
	\caption{Default solver setup.}
	\label{tab:training_solversetup}
\end{table}

\subsection{Training problem set}

We generate difficult random QP problems to train the RL policy $\pi_\zeta$ illustrated in Figure~\ref{fig:framework}.
Each problem of the form \eqref{eq:QP} is defined by matrices $Q, A, G$ and vectors $q, b, d$, with properties such as convexity, sparsity, rank, and conditioning carefully controlled to produce meaningful and diverse problem batches.
A brief overview of these aspects is provided here; implementation details and statistical analyses are available with the accompanying code at \href{https://github.com/JrmBertoncini/DegeneratedQPgenerator.git}{github.com/JrmBertoncini/DegeneratedQPgenerator.git}.

\paragraph{Ill-conditioning}

QP problems become numerically challenging when the linear system \eqref{eq:newton_system} is ill-conditioned, as this directly affects the accuracy of the linear solver in each Newton step.
When the condition number is high, small numerical errors are amplified during factorization, making the Newton directions inaccurate and slowing down convergence, or breaking it altogether.
In our training set, the rank and condition number of matrix $Q$ are targeted to make difficult random QPs, as follows.

\paragraph{$Q$ matrix generator}
Matrix $Q$ is constructed as a positive semidefinite matrix with prescribed rank and spectrum.
An orthogonal matrix $V \in \R^{n \times n}$ is obtained via QR factorization of a Gaussian matrix $ N \in \R^{n\times n}$ with, $N = VR$, where $R \in \R^{n\times n}$ is upper triangular.
Furthermore, let $\Lambda \in \R^{n \times n}$ denote a diagonal matrix containing the eigenvalues of $Q$, whose construction is specified below.
The Hessian is then defined as
$Q = \mathcal{S}_\rho\left(V \Lambda V^\top\right) \in \R^{n\times n}$,
where $\mathcal{S}_\rho$ denotes an operator that enforces symmetry and a density level $\rho$ (so that $Q$ has approximately $\rho n^2$ nonzero components).

A rank $h \le n$ is sampled, and the nonzero eigenvalues are partitioned into $n_L = \lfloor h/2 \rfloor$ large and $n_S = h - n_L$ small values.
For measuring the spectral properties of the cost matrix, we consider
the effective condition number 
\[
    \kappa_e(Q) \coloneqq \frac{\lambda_{\max}(Q)}{\lambda_{\min,+}(Q)},
\]
where $\lambda_{\min,+}$ extracts the smallest positive eigenvalue of a matrix.
Given a target effective condition number $\kappa_e^{\text{target}} \ge 1$, we choose $1< c_0 \ll \kappa_e^{\text{target}}$, draw the eigenvalues as
\[
	\lambda_i \sim \mathcal{U}(1,c_0), \quad i = 1, \dots, n_S,  \qquad \lambda_j \sim \mathcal{U}(0.1, 1) \cdot \kappa_e^{\text{target}}, \quad j = n_S +1, \dots, h,
\]
and assemble the diagonal matrix $\Lambda = \diag(\lambda_1,\dots,\lambda_h,0,\dots,0)$.
This spectral construction induces a spread between the smallest positive and largest eigenvalues, yielding a controlled ill-conditioning. 
This results in a sparse positive semidefinite matrix with $\rank(Q) \leq h < n$ and effective condition number $\kappa_e(Q) \approx \kappa_e^{\text{target}}$.

\paragraph{Constraints linear dependencies}

Numerical challenges arise for Newton-type solvers when dealing with linearly dependent constraints in \eqref{eq:QP}, because, with $\delta_y=\delta_z=0$, linear systems \eqref{eq:newton_system} become singular.
In practice, we generate $m$ equality and $p$ inequality constraints and then a fraction of constraints, $\frac{m}{2}$ and $\frac{p}{2}$, is duplicated and thus linearly dependent, increasing the total number of constraints to $\frac{3m}{2}$ and $\frac{3p}{2}$.
In addition, we consider nearly linearly dependent constraints: some duplicated constraints are perturbed by a very small quantity.
These constraints remain effectively active while introducing near-degeneracy, which is known to be challenging for numerical solvers.

\paragraph{Initial guess: residual control and resampling}

It is essential that the training dataset covers a broad range of the state space to ensure that Algorithm~\ref{alg:rpdip} encounters diverse scenarios during convergence. Coverage can be improved by strategically selecting initial points while keeping generation costs low. If initial points are consistently chosen near feasibility and stationarity, the residuals $(r_{\text{primal}}, r_{\text{dual}})$ will exhibit limited variation, biasing the training distribution toward nearly optimal states.
To avoid this, candidate initial solutions $(x, y, z)$ are selected such that the initial residuals satisfy $\log r_{\text{prim}}^{\text{target}}, \log r_{\text{dual}}^{\text{target}}  \sim \mathcal{U}(-2,2)$ and are found by solving a least-squares problem.

The variable $x\in\R^n$ is randomly generated, and the vectors $b$ and $d$ are constructed, aiming at the target primal residual, as
\begin{equation}\label{eq:primal_target}
    b= Ax - r_{\text{prim}}^{\text{target}} 1_m
    \quad\text{and}\quad
    d = Gx - r_{\text{prim}}^{\text{target}} 1_p + \mu,
\end{equation}
where $\mu \in \R^p$ is a random positive value vector ensuring the inequality direction, and $1_m \in \R^m$ denotes the vector of ones of dimension $m$.
Then, the target dual residual $r_{\rm dual}^{\rm target}$ can be obtained by solving the least-squares problem
\begin{equation}\label{eq:ls_problem_residual}
    \argmin_{y,z} \left( \left\| Qx+q+A^\top y + G^\top z \right\|_2 - r_{\rm dual}^{\rm target}\right) ^2
    \quad
    \stt\quad
    z \geq 0
\end{equation}
with respect to the multipliers $y\in\R^m$ and $z\in\R^p$.
In practice, solving the constrained problem \eqref{eq:ls_problem_residual} is computationally expensive for large-scale dataset generation.
Therefore, the unconstrained counterpart is solved and the inequality multipliers $z$ are then projected onto the nonnegative orthant.

To make the training even lighter, each problem instance is reused with different initial guesses. These additional initializations are sampled randomly around the targeted residual configurations, allowing multiple convergence trajectories to be generated from a single QP instance.

\paragraph{Problem statistics}
The training and validation datasets are summarized in Table~\ref{tab:results_sizes}, which reports the ranges and distributions of the problem dimensions, the controlled effective condition number, and the sparsity pattern of $Q$. These settings are used to generate random QPs for both training and validation.
The training dataset is created with problem sizes ranging from $n_{\min}=5$ to $n_{\max}=50$, in order to generate relatively small and fast to construct QPs.

We introduce in Table~\ref{tab:results_sizes} several regimes of random problems; Standard, LP-like, High-$p$, and Equality-heavy, each designed to stress the solver in different ways.
The Standard regime provides a balanced mix of equality and inequality constraints, serving as a baseline with moderate conditioning.
The LP-like regime, characterized by a nearly vanishing Hessian ($Q \approx 0$) and few or no equality constraints, leads to almost degenerate linear systems typical of linear programs.
The High-$p$ regime introduces a large number of inequality constraints, increasing the size of the active set and the Equality-heavy regime leads to large saddle-point systems that (before regularization) are inherently indefinite.

\begin{table}[tbh]
	\centering%
	\begin{tabular}{l|c}
		\hline
		Setup & $n, \text{ with } \quad n_{\min}=5,\; n_{\max}=50 $ \\
		\hline
		Training         & $n \sim \mathcal{U}(n_{\min}, n_{\max})$ \\
		Validation       & $n \sim \mathcal{U}(n_{\min}, n_{\max})$ \\
		Scale $10\times$ & $n \sim \mathcal{U}(10 \cdot n_{\min}, 10 \cdot n_{\max})$ \\
        \hline
        Density $(Q,A,G)$ & $\rho \sim \mathcal{U}(0.15,0.45)$ \\
        Condition number & $\log\kappa_e^{\text{target}} \sim \mathcal{U}([1,4] \cup [12,22]) $ \\
		\hline
	\end{tabular}
    
	\vspace{0.3cm}
	
	\begin{tabular}{l|c|c}
		\hline
		Regime & $m$ & $p$ \\
		\hline
		Standard 
		& $ \frac{n}{\mathcal{U}(1,3)} $
		& $ n \cdot \mathcal{U}(1,3) $ \\
		LP-like 
		& $\{0,1\}$ 
		& $ n \cdot \mathcal{U}(1,4) $ \\
		High-$p$ 
		& $\{0,1,2\}$ 
		& $n \cdot \mathcal{U}(3,33) $ \\
		Equality-heavy 
		& $ n \cdot \mathcal{U}(0.4,0.9) $ 
		& $0$ \\
		\hline
	\end{tabular}
	\caption{Problem generation setup for training and validation QPs.}%
	\label{tab:results_sizes}%
\end{table}

Figure~\ref{fig:problem_size_sparsity} presents the actual problem sizes and sparsity statistics for the generated datasets, including the number of variables, constraints, and the distribution of nonzeros per row in the constraint matrices.
For reference, the figure also includes a comparison with the problem sizes and sparsity patterns from the Maros--Mészáros benchmark set.
Importantly, these problem-dependent quantities are not provided to the RL controller, which only observes the state vector $\sigma$ composed of convergence and barrier-related information.
This deliberate design choice keeps the policy lightweight and problem-agnostic by relying exclusively on solver convergence indicators such as residual scalars.

The training random problems exhibit a broad range of nonzero values in $Q$ and $[A^\top G^\top]$, which are further scaled by a factor of $100\times$ in the scale $10\times$ validation set.
In terms of sparsity, these generated problems maintain a consistent pattern across the dataset.
In contrast, the Maros--Mészáros benchmark test set displays a wider diversity in both the number of nonzeros and the sparsity patterns, making it challenging to solve as a batch using a single solver configuration.

\begin{figure}[tbh]
    \centering%
   \includegraphics[width=0.8\columnwidth]{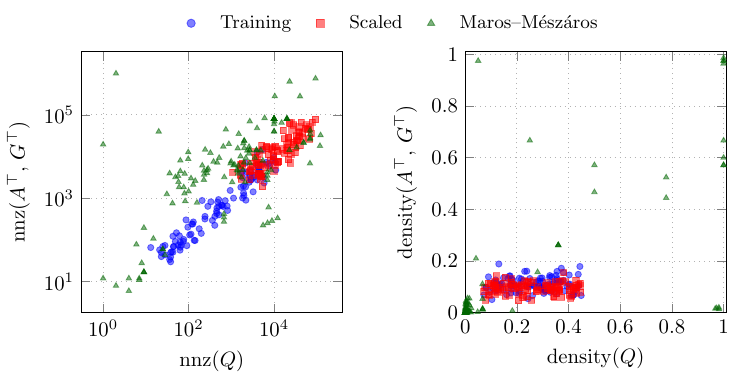}%
    \caption{Distribution of problem size (in terms of number of nonzero entries) and density for randomized training, validation, and scale $10\times$ problems (sample of $10^3$ instances) and for Maros-Mészáros problems (100 smallest instances).}%
    \label{fig:problem_size_sparsity}%
\end{figure}

\begin{figure}[tbh]
	\centering%
    \includegraphics[width=0.8\columnwidth]{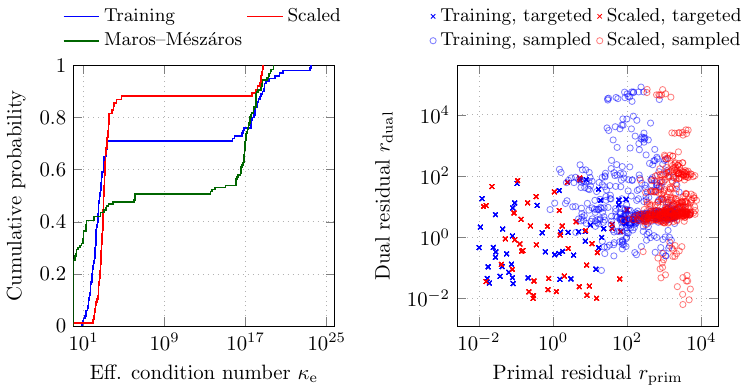}%
    \caption{Left: Empirical cumulative distribution functions (ECDFs) of the effective condition number $\kappa_e$ for the training set, scaled validation set, and Maros--Mészáros benchmark problems. Right: Distributions of the initial primal and dual residuals over a sample of $10^3$ problem instances.}
    \label{fig:problem_kappae_initresiduals}%
\end{figure}

The distribution of effective condition numbers across the training set is shown in Figure~\ref{fig:problem_kappae_initresiduals} (left panel). While a target distribution is prescribed during problem generation, the resulting effective condition numbers differ due to the sparsification operator. Nevertheless, this discrepancy introduces additional variability, which is desirable as it increases the diversity of the training set. The right panel of Figure~\ref{fig:problem_kappae_initresiduals} illustrates the distribution of initial primal and dual residuals across the considered problem classes. The targeted residual values remain within the prescribed range $[10^{-2},10^{2}]$, following a log-uniform sampling strategy. In contrast, the resampled problems tend to exhibit larger residuals overall, with even higher values observed for the scaled validation batch.

\section{Numerical investigations}\label{sec:num_results}

We report the outcomes of the training process and corresponding validation, followed by an evaluation on larger randomly generated problems to assess scalability. Finally, we benchmark the approach on the Maros-Mészáros test set to examine its zero-shot generalization on challenging real-world problem instances. The analysis is based on performance and data profiles, consistently comparing the RL-enhanced solver with the typical fixed-$\alpha$ strategy.
The robustness and acceleration of the different methods will be discussed in the analysis.

\paragraph{Profiles}
For $P$ the set of problems and $S$ the set of solvers, let $t_{s,p}$ denote the (wall-clock) runtime required by solver $s\in S$ to solve instance $p\in P$ (for a given tolerance $\epsilon>0$).
We will monitor the computational performance of different solver settings and graphically summarize our results by means of two profiles.
When solver $s\in S$ fails to solve $p\in P$, we set $t_{s,p}=\infty$.
\begin{itemize}
    \item A \emph{data profile} is the graph of the cumulative distribution function $f_s \colon [0,\infty) \to [0,1]$ of the runtime, namely $f_s(t) \coloneqq \frac{|\{ p\in P \,|\, t_{s,p} \leq t \}|}{|P|}$, where $|P|$ denotes the cardinality of set $P$.
    As such, each data profile reports the fraction of problems $f_s(t)$ that can be solved by $s$ with a computational budget $t$, and therefore it is independent of the other solvers.
    %%%
    \item A \emph{performance profile} is the graph of the cumulative distribution function of the performance ratio $r_{s,p} \coloneqq \frac{t_{s,p}}{\min_{s^\prime \in S} t_{s^\prime,p}}$.
    As such, it displays the fraction of problems solved within a time factor compared to the fastest solver of each problem.
    Therefore, the performance profile of a solver also depends on the performance of all other solvers in $S$.
\end{itemize}

\subsection{Training results}

The overall training spans under 8 hours and includes 800k randomly generated QP problems.
The training progress can be observed in Fig.~\ref{fig:trainingstats}, where the light noisy curves show the raw per-episode values whereas the bold line is a moving-window average.
As metrics to monitor the training, we observe the value function loss, the policy loss, the reward and the entropy in each episode. 
The loss value tends to show the training behavior of the policy and the value function.
The reward represents the objective of the process and therefore it needs to be maximized.
The entropy is monitored, since it represents the exploration and uncertainty of the policy during the training.

In Fig.~\ref{fig:trainingstats_vf} it can be observed how the loss of the critic evolves.
In first 600k steps a significant reduction of the loss can be noticed, which indicates that the approximation of the value function improves.
In the last 200k steps, the PPO algorithm fine-tunes this approximation.
In addition to the loss function, Fig.~\ref{fig:trainingstats_pl} shows the policy loss, which is an indicator of the improvement of the policy.
The mostly negative values show that the policy improves, since advantageous actions, which means that they have a positive advantage function approximation, are reinforced, and disadvantageous actions are weakened.

Fig.~\ref{fig:trainingstats_rew} shows the evolution of the episode reward throughout training.
After an initial exploration dip at $\approx$25k steps, the policy recovers and improves monotonically.
The variance of the raw signal reflects the probabilistic behavior of the environment while the smoothed trend shows no plateau within the reported horizon.

Furthermore, we plot the entropy during training in Fig.~\ref{fig:trainingstats_ent}.
Entropy can be used to determine how uncertain or exploratory the training is. 
We see that the entropy increases first, which is a typical behavior, since exploration is needed at the beginning of the training.
Afterwards, the entropy stagnate to a certain level.
Although this behavior is common for RL algorithms in general, it is noticeable that the final level of the entropy is relatively high. 
This means that the final policy is relatively uncertain and is not close to a deterministic policy.
The investigation of whether an enhanced state space or even longer training phases can make the policy more certain is left for future works.
In any case, the plots indicate a successful training, and the influence of the RL policy on the QP solver is discussed in the remainder of this manuscript.

\begin{figure}[tbh]
	\centering%
	\begin{subfigure}{0.47\columnwidth}
		\centering
		\includegraphics[width=\linewidth]{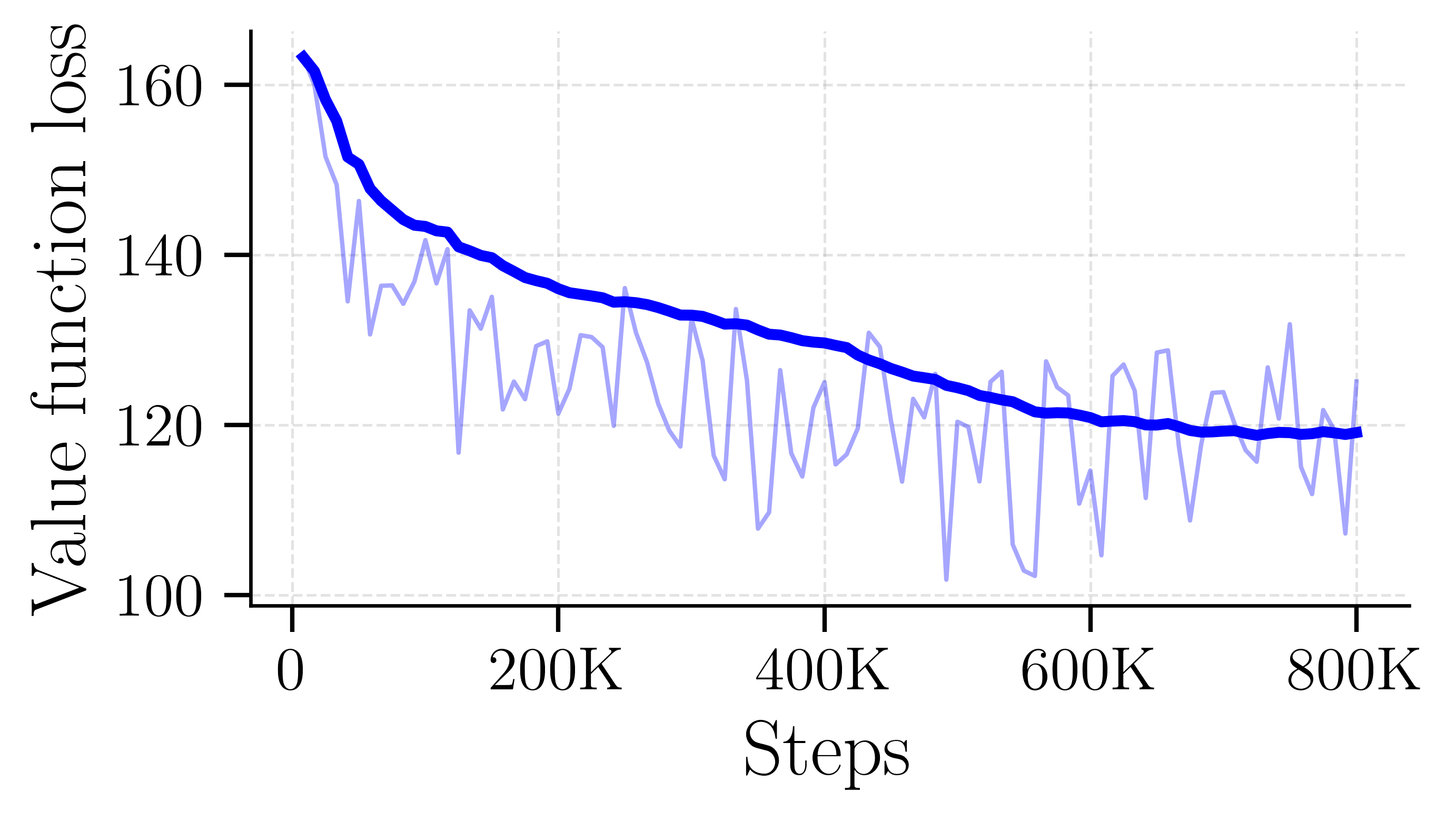}
		\caption{Value function loss}
		\label{fig:trainingstats_vf}
	\end{subfigure}
	\hfill
	\begin{subfigure}{0.47\columnwidth}
		\centering
		\includegraphics[width=\linewidth]{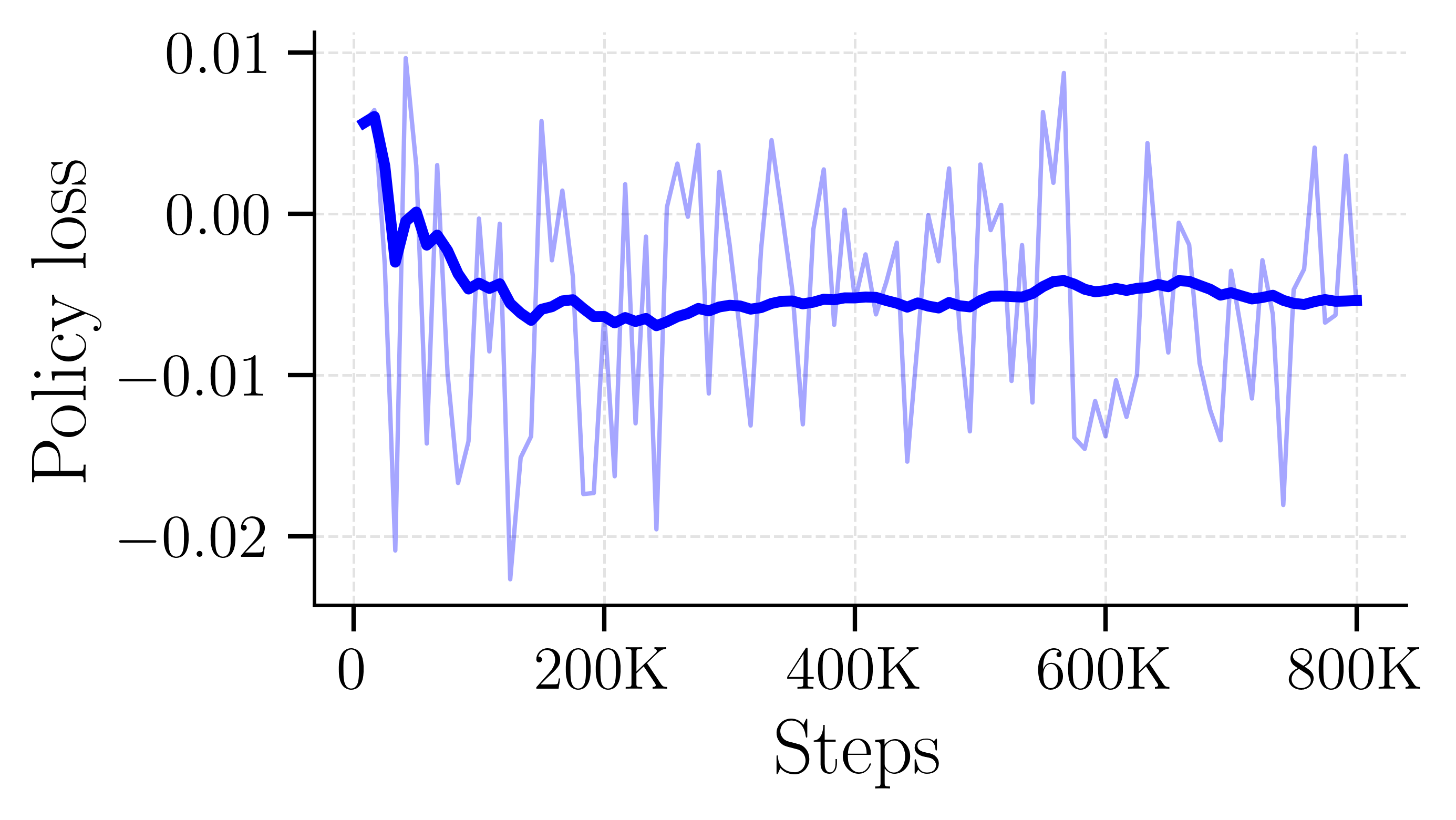}
		\caption{Policy loss}
		\label{fig:trainingstats_pl}
	\end{subfigure}
    \\
	\begin{subfigure}{0.47\columnwidth}
		\centering
		\includegraphics[width=\linewidth]{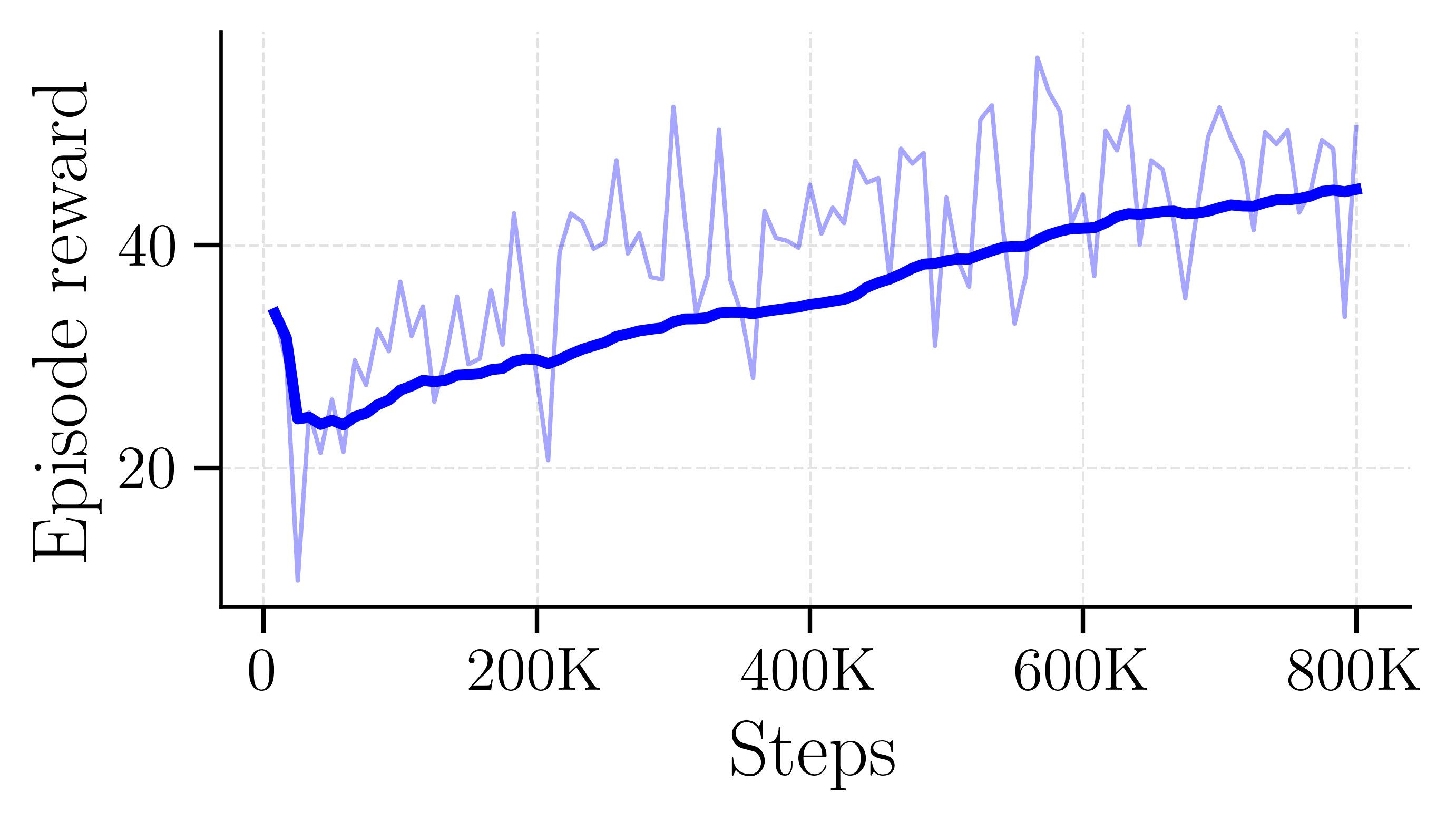}
		\caption{Reward}
			\label{fig:trainingstats_rew}
	\end{subfigure}
     \hfill
	\begin{subfigure}{0.47\columnwidth}
	\centering
	\includegraphics[width=\linewidth]{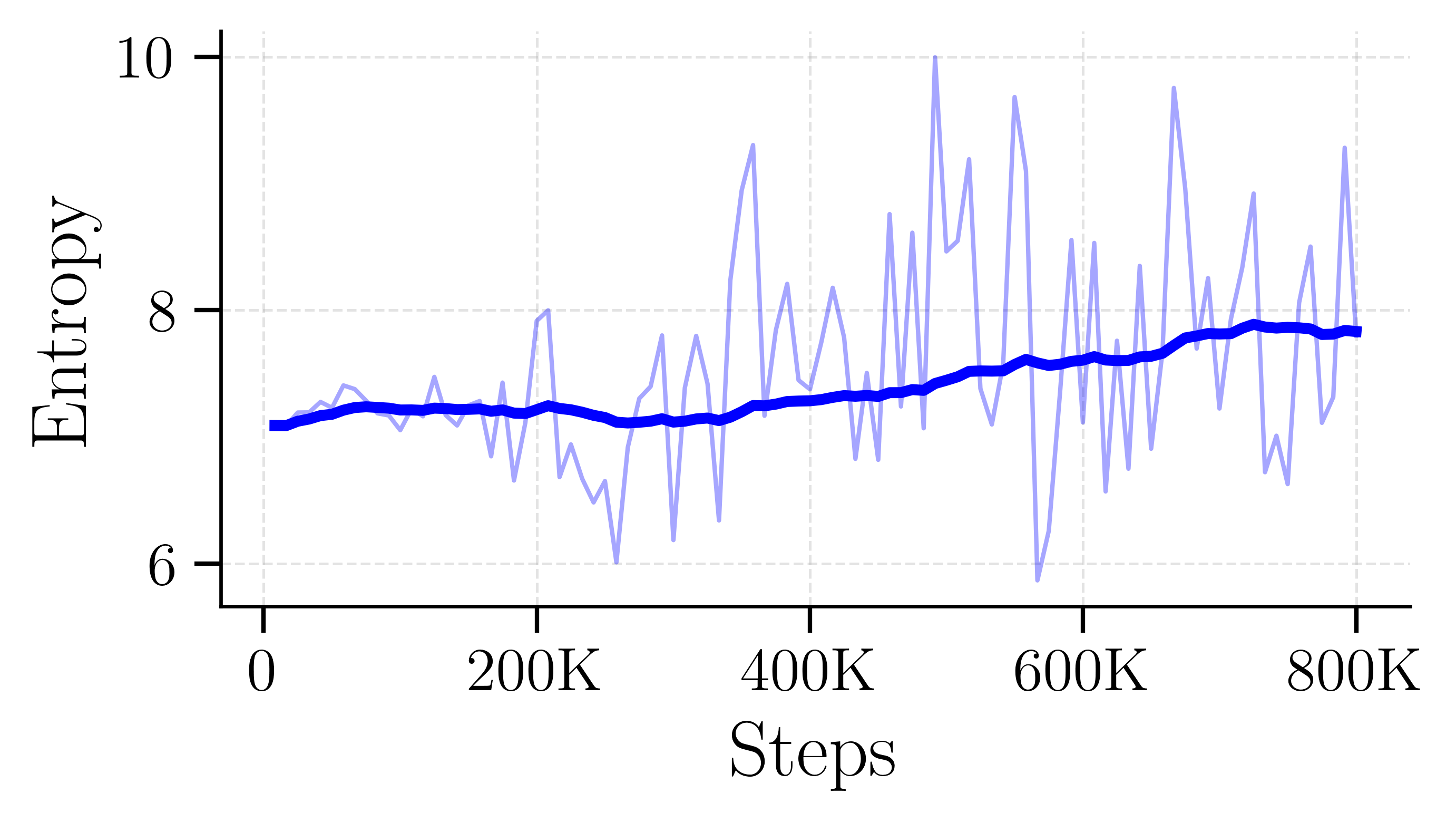}
	\caption{Entropy}
		\label{fig:trainingstats_ent}
    \end{subfigure}
	\caption{Evolution of relevant metrics during the training procedure: raw per-episode values (light blue) and moving-window average (dark blue).}%
	\label{fig:trainingstats}%
\end{figure}

\subsection{Validation}

After the training, the performance of the method is evaluated on a batch of $100$ problems having the same statistical properties as the problems used for the training, (see Table~\ref{tab:results_sizes}).
Figure~\ref{fig:profiles_validation} reports the performance and data profiles
on the validation set, which is drawn from the same distribution as the training
batch and consists of small-scale QPs with $n \in [5, 47]$ (median $n = 17$).

At performance ratio $r = 1$, \rlqpsolver{} solves $\approx 50 \%$ of
problems as the best solver, against $ \approx 35\%$ for the baseline \baseqpsolver{},
a relative gain of $\approx 30\%$ in problems solved at
parity. Both curves saturate at the same ceiling of $87.0\%$ for $r \geq 2$,
confirming that the RL policy does not recover any previously unsolvable problem but consistently reaches the solution faster.

The data profile shows that all solvable problems are handled well under one
second, with both solvers saturating at $t = 0.5$\,s.
The bulk of the solve time is concentrated around $100$\,ms (median runtime
$\approx 65$\,ms for both solvers.
At $t = 50$\,ms, \rlqpsolver{} has already solved $\approx 12\%$ of problems
compared to $\approx 5\%$ for the baseline. The inference overhead of the neural policy (measured in the running time) recorded amounts to $\approx 1$\,ms per problem, which is negligible relative to the solve time.

On problems of this scale, the iteration counts are already low leaving little room for large absolute savings.
A consistent gain of over $30$ percentage points in problems solved at parity
therefore represents a meaningful improvement: the RL policy systematically
selects parameters that reduce the number of innerloop steps even on
small, ill-conditioned instances where a hand-tuned baseline already performs
well.

\begin{figure}[tbh]
	\centering%
	\includegraphics[width=0.8\columnwidth]{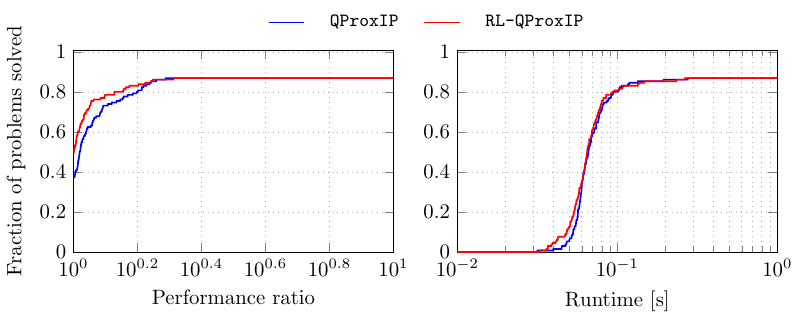}%
	\caption{Performance (left) and data (right) profiles of a validation test-set.}%
	\label{fig:profiles_validation}%
\end{figure}

Figure~\ref{fig:alphamaps} shows a 3D voxel-based state-to-action mapping of the learned policy defined over the observation space.
Each voxel aggregates all samples falling into the corresponding region of the discretized state space, and the displayed value corresponds to the average of the associated action component over these samples.
The marker size is proportional to the number of samples per voxel.
Since \(\varepsilon\) represents the inner-loop tolerance, its value decreases as the solver progresses. 

During convergence, the residuals decrease from large to small values, while $\varepsilon$ decreases. The resulting trajectories cluster around a hyperplane, as illustrated in Figure~\ref{fig:alphamaps}. Along these convergence paths, the learned policy exhibits a clear structure across all action components. In particular, $\alpha_x$, $\alpha_y$, and $\alpha_z$ remain within approximately $[0.05, 0.3]$, while $\alpha_{\nu}$ lies in $[0.2, 0.4]$, and $\alpha_{\varepsilon}$ is consistently higher, in the range $[0.7, 0.95]$. This separation suggests differentiated roles of the control parameters across solver regimes.

\begin{figure}[tbh]
    \centering%
    \includegraphics[width=0.95\columnwidth]{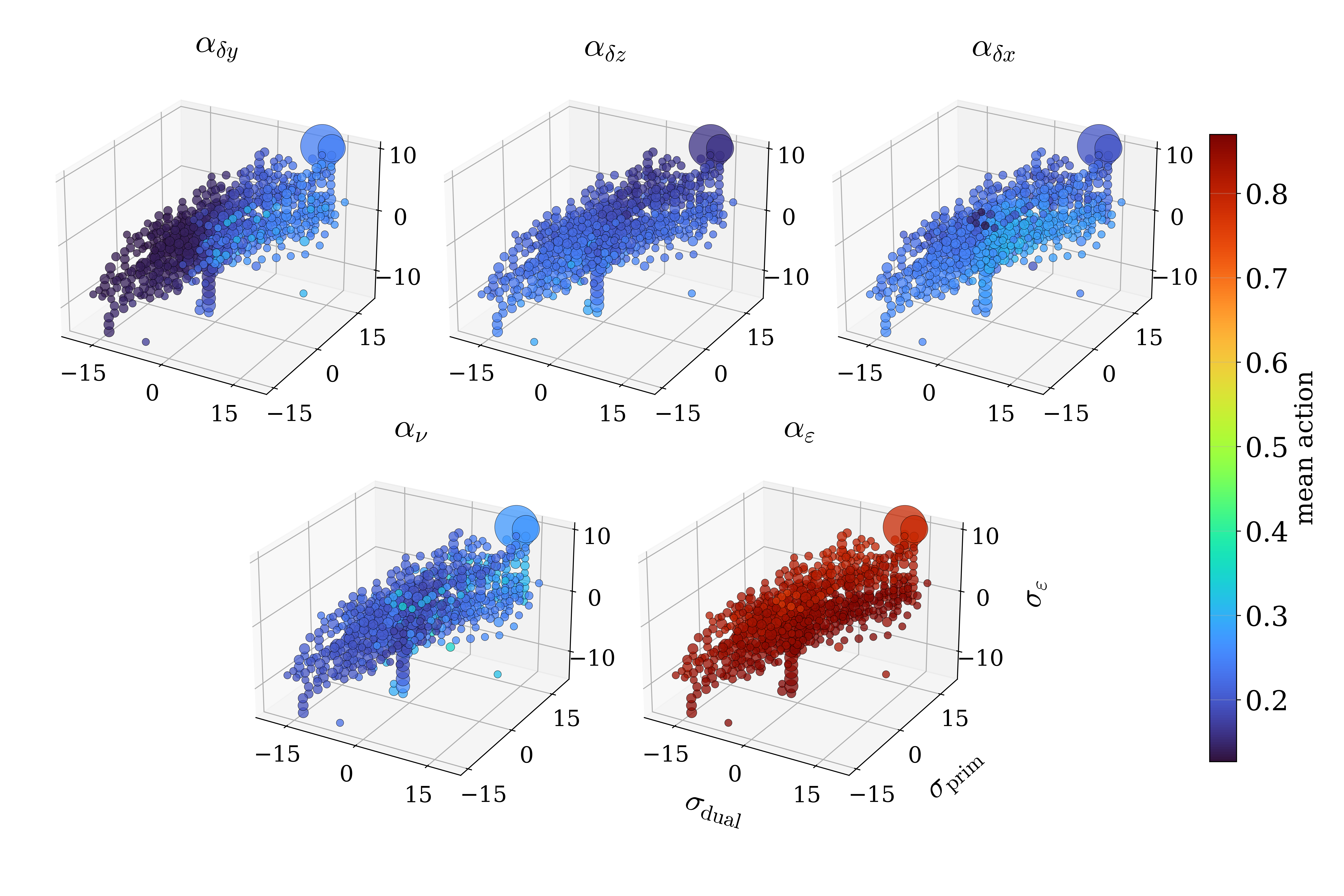}%
    \caption{Mappings $\alpha$ as a function of $\sigma_{\rm prim}$, $\sigma_{\rm dual}$, and $\sigma_{\varepsilon}$ evaluated on the 100 smallest Maros--Mészáros problems. The dependency on $\nu$ is not shown.}
    \label{fig:alphamaps}
\end{figure} 

\subsection{Problem dimensions scaling}

When evaluating the RL-policy on the validation set, the observed improvements over the base solver \baseqpsolver{} (with fixed $\alpha$ factors, $\alpha_x=0.15, \alpha_x=0.2$ and $\alpha_x=0.2$) remain modest.
This is largely due to the relatively small size of the validation problems, which, although ill-conditioned, non-strictly convex, and subject to redundant constraints, still exhibit limited complexity.
However, the choice of random problem sizes during training  (see Figure~\ref{fig:problem_size_sparsity} for a comparison of problem properties) was primarily motivated by the need to keep the computational training cost low enough to complete on a personal computer within a single day; see the system configuration in Table~\ref{tab:PPO_hyperpara}.
A question then arises: How does the policy $\pi_\zeta$, trained at low cost on small QPs, perform on larger problem instances?

Figure~\ref{fig:profiles_scale10} presents the results for a batch of problems whose dimensions and sparsity are scaled by a factor of $10$ relative to the validation set.
To assess the performance of \baseqpsolver{} and \rlqpsolver{}, we compare it against several widely used quadratic programming solvers representing distinct algorithmic paradigms.
\texttt{PIQP} is a modern primal-dual interior-point solver designed specifically for quadratic programs.
It employs a Mehrotra predictor-corrector strategy to achieve high numerical robustness and efficiency, particularly for high termination tolerance, see \cite{schwan2023piqp}.
The \texttt{OSQP} solver implements ADMM for convex quadratic programs.
It is well suited for large, sparse problems, offering predictable performance and simple warm-start capabilities, though it may require more iterations to reach high-accuracy solutions, see \cite{stellato2020osqp}.
\texttt{ProxQP} is a proximal-based quadratic programming solver leveraging iterative refinement and proximal regularization.
It provides robust convergence behavior and predictable numerical performance, particularly suitable for embedded or control-oriented applications, \cite{bambade2025proxqp}. 

We note that \texttt{OSQP} is not included in the profiles because it did not converge to the required tolerance for any of the test problems. This is consistent with its design as a first-order ADMM-based solver, which typically struggles to achieve high accuracy. We further emphasize that no preconditioner is used for any of the solvers.

\begin{figure}[tbh]
	\centering%
	\includegraphics[width=0.8\columnwidth]{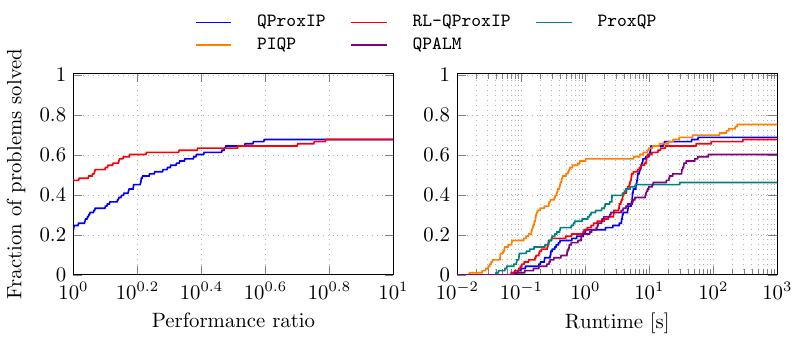}%
	\caption{Scale $10\times$ set: Performance (left) and data (right) profiles comparing benchmarked solvers on a 10$\times$ scaled training set.}
	\label{fig:profiles_scale10}%
\end{figure}

\subsection{Maros--Mészáros problems}

Training was carried out on randomly generated problem instances, but an important question is how well the method generalizes to established benchmarks. To this end, we evaluate on the Maros–Mészáros test set \cite{maros1999repository}, a widely used benchmark collection of challenging quadratic programs. This suite contains many large-scale and ill-conditioned problems and is considered a standard for assessing the robustness and efficiency of convex QP solvers. For evaluation, we consider the problem in a zero-shot generalization setting and focus on the 100 smallest instances (out of 138), where problem size is determined by folder size.
A comparison of these instances with the randomly generated problems used for training is shown in Figure~\ref{fig:problem_size_sparsity}.
As before, the Maros--Mészáros set is used in its raw form, meaning that inequality constraints with ``infinite'' bounds (represented as values of $10^{20}$ in OSQP \cite{stellato2020osqp}) are retained in the problem data \eqref{eq:QP}.
To accommodate the increased difficulty of these benchmarks, both the maximum number of inner and outer iterations are increased to $50$.

Fig.~\ref{fig:profiles_validation_MM} reports the performance and data profiles for \baseqpsolver{} and \rlqpsolver{}.
\rlqpsolver{} is the fastest solver on $\approx 74\%$ of instances (performance ratio $r = 1$), while \baseqpsolver{} requires $r \approx 4$ to reach the same fraction.
The curves cross at $r \approx 10$: \baseqpsolver{} solves $100\%$ of instances while \rlqpsolver{} saturates at $\approx 92\%$.
In absolute runtime, \rlqpsolver{} solves $50\%$ of instances within $t \approx 40\,\mathrm{ms}$, roughly $2.3\times$ faster than \baseqpsolver{}.

\begin{figure}[tbh]
	\centering%
	\includegraphics[width=0.8\columnwidth]{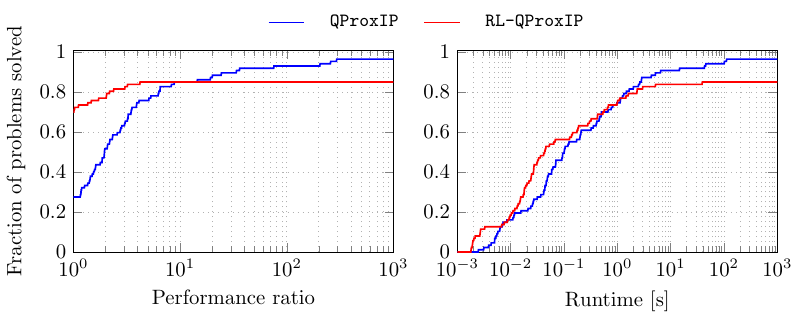}%
	\caption{Performance (left) and data (right) profiles on the Maros--Mészáros problems.}%
	\label{fig:profiles_validation_MM}%
\end{figure}

To further analyze the effect of the RL policy, we restrict the benchmark to the subset of problems for which \rlqpsolver{} terminates within one second.
This threshold is approximately $10^3$ times longer than the time required to solve the small random problems during the training phase.
Thus, it provides a genuine out-of-distribution evaluation.

\begin{figure}[tbh]
	\centering
		\centering%
		\includegraphics[width=0.8\columnwidth]{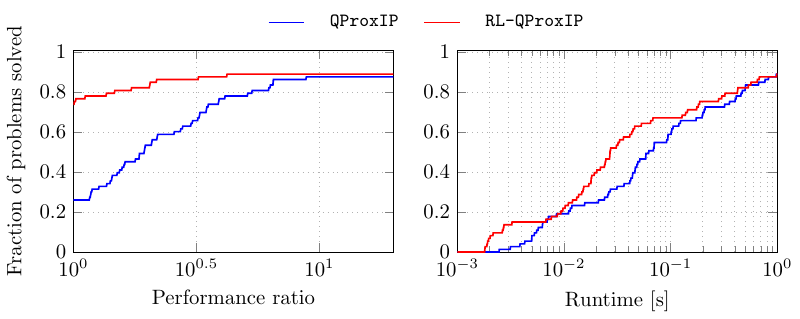}%
        \caption{Performance (left) and data (right) profiles on the subset of 73 Maros--Mészáros problems solved by at least one of the two solvers (\baseqpsolver{} or \rlqpsolver{}) within the one-second time budget.}%
	\label{fig:profiles_validation_MM_1sec}%
\end{figure}

Figure~\ref{fig:profiles_validation_MM_1sec} presents the performance and data profiles on the subset of 73 Maros--Mészáros problems that at least one of the two solvers (\baseqpsolver{} or \rlqpsolver{}) solves within a one-second time budget. We exclude only the instances that neither solver can handle under the one second time limit. The performance profile indicates that \rlqpsolver{} is the fastest solver on 74\% of the problems (performance ratio $r = 1$), compared to 26\% for \baseqpsolver{}. In absolute terms, \rlqpsolver{} solved 66 of the 73 instances, while \baseqpsolver{} solved 65. These results demonstrate \rlqpsolver{} generalizes effectively beyond its training distribution, delivering consistent speedups even on significantly larger and more challenging problems.

\begin{table}[tbh]
    \centering
    \begin{tabular}{lrrr}
        \hline
        \textbf{Problem} & \baseqpsolver{} [s] & \rlqpsolver{} [s] & \textbf{Speedup} \\
        \hline
        PRIMALC8   & 6.433  & 0.0219 & 294 \\
        PRIMALC2   & 5.214  & 0.0203 & 257 \\
        PRIMALC5   & 2.005  & 0.0098 & 205 \\
        PRIMALC1   & 1.146  & 0.0151 & 76 \\
        DUALC8     & 0.9782 & 0.0271 & 36 \\
        QSHIP08S   & 101.1  & 2.999  & 34 \\
        QSHIP04S   & 1.666  & 0.0683 & 24 \\
        QSHIP04L   & 2.867  & 0.1436 & 20 \\
        QETAMACR   & 45.43  & 2.371  & 19 \\
        CVXQP2\_M  & 14.71  & 1.013  & 15 \\
        \hline
    \end{tabular}
    \caption{Comparison on problems where the RL-enhanced solver \rlqpsolver{} is at least $10\times$ faster than the base solver \baseqpsolver{}.}
    \label{tab:extreme_speedups}
\end{table}

Table \ref{tab:extreme_speedups} shows that \rlqpsolver{} achieves substantial acceleration over \baseqpsolver{} on a small subset of problems, with speedups reaching almost $300 \times$ in the most extreme cases. Notably, a clear pattern emerges across specific problem classes—particularly \textsc{PRIMALCx} and \textsc{QSHIP0xx} instances, where \rlqpsolver{} consistently delivers significantly reduced solve times compared to the baseline solver.

\section{Final remarks and outlook}\label{sec:conclusions}

We presented a method based on reinforcement learning to speed up and harden the convergence of a regularized interior-point solver for convex quadratic programming.
We used RL to learn a policy aimed at adapting the solver's hyperparameters to reduce the number of iterations needed to obtain a solution.
For the RL training, we generated a framework to design small random QPs, which serves as good training problems.
Then the parameterized RL policy, which is constructed as an artificial neural network, is trained on these QPs.
In experiments, the results suggest that a single policy can improve convergence rates for a broad class of problems.
In particular, the policy exhibits an input-output behavior that is interpretable, based on the underlying QP method, and it can successfully be used for problems with different sizes, scales and structures. 
Compared to previous approaches, like RLQP \cite[\S 6]{ichnowski2021rlqp}, our lightweight training achieves similar validation performance with a fraction of time and memory (7 hours on a laptop rather than days on a high-end computer).
For QPs that are solved in only a few iterations and thus require limited hyperparameter adaptation, a learned policy may provide only marginal additional benefit.
In these cases, the inference overhead of evaluating the RL policy can partially offset the gains from faster convergence. Nevertheless, already for problems with sub-second solve times, \rlqpsolver{} is competitive with state-of-the-art solvers while achieving speed-ups at a comparable level of robustness.

\medskip

Future work will investigate the policy performance to other QP problem classes, such as large-scale sparse problems arising in power systems, optimal control problems showing particular three-banded structures or financial optimization with highly volatile market data.
Moreover, for even more adaptive RL approaches, the policy could continue learning and improving during solver deployment.
Moreover, we will assess whether the presented RL policy, trained within a primal–dual interior-point framework, transfers to other second-order optimization schemes such as primal interior-point methods.
More broadly, this raises the question of generalization across Newton-type algorithms, such as between augmented Lagrangian and interior-point methods.
Such an analysis would help clarify to what extent the learned RL strategy captures algorithm-agnostic structure of second-order optimization dynamics, as opposed to being specific to a particular solver architecture.

%---------- References
\phantomsection
\addcontentsline{toc}{section}{References}%
\bibliographystyle{habbrv}
\bibliography{biblio}

\end{document}